\documentclass{article}

\usepackage{geometry}
\usepackage{lmodern} 
\usepackage{wrapfig}
\usepackage{colortbl}
\usepackage{amsmath}
\usepackage{amssymb}
\usepackage{mathrsfs}
\usepackage{amsthm}
\usepackage{makeidx}

\usepackage{multirow}
\usepackage{caption}
\usepackage{arydshln}

\newtheorem{thm}{Theorem}[section]

\usepackage{graphicx}

\usepackage{float}

\floatstyle{ruled}
\newfloat{triangle}{!ht}{tar}
\floatname{triangle}{arithmetic triangle}

\usepackage{ulem}

\usepackage{url}
\usepackage{lipsum}

\newcommand\blfootnote[1]{%
  \begingroup
  \renewcommand\thefootnote{}\footnote{#1}%
  \addtocounter{footnote}{-1}%
  \endgroup
}

\date{} 

\title{The $3x+1$ and $5x+1$ Problems}

\author{Robert Tremblay}

\begin{document}


\maketitle 


\begin{abstract}
We will prove that there are trajectories generated by the function at the origin of the $5x+1$ problem which are divergent. The iterative application of this function on the set of positive integers allows us to determine that more than $17$ $\%$ of all these integers start divergent trajectories. Regarding the $3x+1$ problem, this percentage tends towards zero, suggesting that all positive integers are part of converging trajectories. Despite this appearance, we cannot conclude that all positive integers belong to convergent trajectories. Nevertheless, the results obtained in this paper allow us to follow the evolution of the distribution of trajectories and to understand why the values of positive integers ending in the integer $1$ are more and more large. 

\end{abstract}

\section{Introduction}
 
Let $q$ a integer, positive, negative or zero. We  define the function $f(2q) = q/2$ applied to even integers $2q$ and the functions $f_3(2q+1) = 3q+2$ ($3x+1$ problem) or $f_5(2q+1) = 5q+3$ ($5x+1$ problem) those applied to odd integers $2q+1$. The iterative application of these functions on the set of integers generates sequences of integers called trajectories. By using this process of iterations on different integers one quickly observes a chaotic behavior. In fact, there does not seem to be any correlation between the different types of trajectories and the starting integer of the sequence. The same is true for the values reached.

Before developing the method that we will use in this paper to deal with these problems, remember that the $3x+1$ conjecture states that if we start with from any positive integer and we iteratively apply the functions $f$ and $f_3$ above, we eventually end up with $1$. If the conjecture is false, it means that there are positive integers $n \ge 2$ which never satisfy the condition $f^{(k)}(n)< n$ when the number of iterations $k$ approaches infinity.

First, we develop a method to deal with all the integers as large as we want. This is made possible because there is a periodicity in the distribution of trajectories.

Subsequently, the essence of our approach is based on the notion of distribution function of integers starting the different trajectories and, the detailed analysis of the properties of this function will naturally lead us to the assertions stated in the abstract.

\section{Functions $T_3$ and $T_5$ and periodicity}

Mappings can be define on integers represented by functions such that each element of the set $\mathbb{Z}$ is connected to a single element of this set. The iterative application of these functions produces a sequence of integers called trajectories. 

Let

\begin{equation*}
(n, f(n),f^{(2)}(n),f^{(3)}(n), \cdots, f^{(i)}(n), \cdots) ,
\end{equation*}

 with $f^{(i+1)}(n) = f\{f^{(i)}(n)\}$, $i =0, 1, 2, 3, \cdots$ and $f^{(0)}(n) = n$, a trajectory generates by a function f on an integer $n$.

A sequence of integers forms a loop when there exists a number of iterations $k \ge 1$ such that 

\begin{equation}
f^{(k)}(n) = n. 
\label{condition_cycle}
\end{equation}

If all integers in the sequence are different two by two, we have by definition a cycle of length $p = k$, so the trajectory $(n, f(n),f^{(2)}(n),f^{(3)}(n), \cdots, f^{(k-1)}(n))$. Generally, we note the trajectory characterizing a cycle starting with the smallest integer. 

There are a multitude of functions that have these properties. The function $g(n)$ giving rise to the original Collatz problem and the $3x+1$ function $T_3(n)$~\cite{lagarias}, the $5x+1$ function $T_5(n)$ and the accelerated $3x+1$ function~\cite{kontorovich_lagarias}, are some examples. Except for the last function, the others come from a group called Generalized $3x+1$ Mappings~\cite{matthews}.

The two functions dealt with in this paper are defined by 

\begin{equation}
T_{m_i}(n) = \left\lbrace
\begin{array}{ll}
\frac{n}{2} & \mbox{, if $n\equiv0\pmod{2}$},\\
\\
\frac{m_in+1}{2} & \mbox{, if $n\equiv1\pmod{2}$}\\
\end{array}
\right.
\label{mx+1_function}
\end{equation}

with $m_3 = 3$ for the $3x+1$ problem and , $m_5 = 5$ for the $5x+1$ problem.

The general expression giving the result of $k$ iterations of the function $T_{m_i}$, which we will simply call $T$, on an integer $n$ is 

\begin{equation}
	T^{(k)}(n) = \lambda_{k_1,k_2}n + \rho_k(n),  
\label{g_k}
\end{equation}

where

\begin{equation}
	\lambda_{k_1,k_2} = \left(\frac{1}{2}\right)^{k_1}\left(\frac{m_i}{2}\right)^{k_2} 
\label{lambda_k1_k2}
\end{equation}

and

\begin{equation}
	k = k_1 + k_2,
\label{k1_plus_k2}
\end{equation}

with $k_1$ the number of transformations of the form $n/2$ and $k_2$, transformations of the form $(m_in + 1)/2$.

Unlike parameter $\lambda_{k_1,k_2}$, $\rho_k(n)$ depend on the order of application of the transformations.

Let $n$ and $T^{(k)}$ be replaced by the variables $x$ and $y$,

\begin{equation*}
	2^k\rho_k(n) = 2^ky - 1^{k_1}m_i^{k_2}x.  
\label{eq_diophantine}
\end{equation*}

In this form we have a Diophantine equation of first degree at two unknowns, 

\begin{equation}
	c = by - ax,  
\label{eq_diophantine_general}
\end{equation}

where 

\begin{equation}
	a = 1^{k_1}m_i^{k_2}, \phantom{12} b = 2^k \phantom{12} and \phantom{12} c = 2^k\rho_k(n).  
\label{par_a_b_c}
\end{equation}

Depending on the new parameters $a$ and $b$ the parameter $\lambda$ (equation \ref{lambda_k1_k2}) becomes 

\begin{equation}
	\lambda_{a,b} = \left(\frac{a}{b}\right). 
\label{lambda_a_b}
\end{equation}

From a well-known result of Diophantine equations theory we have the theorem

\begin{thm}
Let the Diophantine equation $c = by - ax$ of first degree at two unknowns. If the coefficients $a$ and $b$ of $x$ and $y$ are prime to one another (if they have no divisor other than $1$ and $-1$ in common), this equation admits a infinity of solutions to integer values. If $(x_0, y_0)$ is a specific solution, the general solution will be $(x = x_0 + bq, y = y_0 + aq)$, where $q$ is any integer, positive, negative or zero.
	\label{equationDiophantine}
\end{thm}

\textsl{Proof}

References : Bordell\`es~\cite{bordelles}. $\blacksquare$

We may to assign to every integer of a trajectory generates by the function $T(n)$ a number $t_j = 0$ if $T^{(j)}(n)$ is even, and $t_j = 1$ if it is odd. Then, the iterative application of the function $T$ to an integer $n$ give a dyadic sequence $w_l$ of $1$ and $0$  


\begin{equation*}
w_l = (t_0, t_1, t_2, t_3, \cdots, t_j, \cdots, t_{l-1}), \phantom{1,2} with \phantom{1,2} l \ge 1.
\end{equation*} 


For a given length $l$ there are $2^l$ different dyadic sequences $w_l$ of $0$ and $1$. 

The representation of the trajectories in terms of $t_{j}$ leads to an important theorem which makes it possible to bring out an intrinsic property, namely the \textsl{periodicity}. This property has already been observed by Terras~\cite{terras} and Everett~\cite{everett} concerning the process of iterations of the function $T_3(n)$ generating the problem $3x+1$, and appears in a theorem which they have demonstrated by induction. We will prove it differently, using the previous theorem.


\begin{thm}
	All dyadic sequences $w_l$ of length $l = k \ge 1$ generated by any $2^l$ consecutive integers are different and are repeated periodically.
	\label{periodicity}
\end{thm}

\textsl{Proof}

Let $k = l \ge 1$ the number of iterations applied to a given integer $n$. The trajectories of length $L = l + 1$

\begin{flushleft}
	$\phantom{1,2,3,4} (T^{(0)}(n),T^{(1)}(n))$ 
	
	$\phantom{1,2,3,4} (T^{(0)}(n),T^{(1)}(n),T^{(2)}(n))$ 
	
	$\phantom{1,2,3,4} \cdots$ 
	
	$\phantom{1,2,3,4} (T^{(0)}(n),T^{(1)}(n)), \cdots, T^{(k)}(n))$
\label{trajectories}
\end{flushleft}

correspond respectively to the dyadic sequences of length $l \ge 1$ 

\begin{flushleft}
	$\phantom{1,2,3,4} w_1 = (t_0)$ 
	
	$\phantom{1,2,3,4} w_2 = (t_0, t_1)$ 
	
	$\phantom{1,2,3,4} \cdots$
	
	$\phantom{1,2,3,4} w_{l=k} = (t_0, t_1, \cdots, t_{k-1})$ .
\label{sequences}
\end{flushleft}

For a given number $l$ we have $2^l$ different dyadic sequences $w_l$ possible. 

According to theorem~\ref{equationDiophantine}, each of the $2^l$ dyadic sequences will be performed for $k = l$. Indeed, the $0$ and the $1$ of these sequences correspond to the operations on the even and odd integers. We build $2^k$ different Diophantine equations characterized by $2^k$ different combinations of the parameters $a$, $b$ and $c$, whose solutions will be given by $(x = x_0 + 2^kq, y = y_0 + m_i^{k_2}q)$. Therefore, all the integers $x_0 + 2^kq$ starting a trajectory of length $k + 1$ correspond to the same sequence $w_k$. In a sequence of $2^k$ consecutive integers, each integer must start a different sequence $w_k$, otherwise the $2^k$ different dyadic sequences will not be performed. $\blacksquare$

Consequently, the periodicity property allows us to distribute all the integers starting the different trajectories in slices of $2^k$ consecutive integers. Like all dyadic sequences are realized, the distribution is binomial. If $k$ is the total number of iterations composed of $k_1$ operations on even integers and $k_2$ operations on odd integers, then the number of dyadic sequences is obtained by the binomial coefficients (BC), so

\begin{equation}
		BC = \left( \begin{array}{c} k \\ k_2 \end{array}\right) = \left( \begin{array}{c} k \\ k_1 \end{array}\right) = \frac{k!}{k_1!k_2!}.
\label{BC}
\end{equation}  

These coefficients can be calculated otherwise, either by using a property inherent in Pascal's triangles. By numbering each of the rows $i = k_2$ and each of the columns $j =k$, we can write $n(i,j) = BC$, so the different elements of the table~\ref{PascalTriangle_BC} (table
 \ref{PascalTriangle_BC_to_k_EQ_20} in appendix). By the properties of the binomial coefficients we have $n(0,k) = 1$ (the top line) and $n(k,k) = 1$ (the bottom diagonal). In addition, all other elements are the result of sum 

\begin{equation}
	n(0 < i < j, j > 1) = n(i-1,j-1) + n(i,j-1)	
\label{somme}
\end{equation}

which corresponds to a recursive calculation.


\begin{table}[H]
\begin{center}
\begin{tabular}{c|cccccccccccc}
	
	$k_2 \setminus k$  & 0 & 1 & 2 & 3 & 4 & 5  & 6  & 7  & 8  & 9 & 10 &$\cdots$ \\
		\hline
	0 & 1 & 1 & 1 & 1 & 1  & 1  & 1  & 1  & 1  & 1 & 1 & $\cdots$ \\
	1 &   & 1 & 2 & 3 & 4  & 5  & 6  & 7  & 8  & 9 & 10 & $\cdots$ \\
	2 &	  &		& 1 & 3 & 6  & 10 & 15 & 21 & 28 & 36 & 45 & $\cdots$ \\
	3	&		&		&   & 1 & 4  & 10 & 20 & 35 & 56 & 84 & 120 & $\cdots$ \\
	4	&		&		&		&   & 1  & 5  & 15 & 35 & 70 & 126 & 210 & $\cdots$ \\
	5	&		&		&		&   &    & 1  & 6  & 21 & 56 & 126 & 252 & $\cdots$ \\
	6	&		&		&		&   &    &    & 1  & 7  & 28 & 84 & 210 & $\cdots$ \\			
	7	&		&		&		&   &    &    &    & 1  & 8  & 36 & 120 & $\cdots$ \\
	8	&		&		&		&   &    &    &    &    & 1  & 9 & 45 & $\cdots$ \\
	9	&		&		&		&   &    &    &    &    &    & 1 & 10 & $\cdots$ \\
	10	&		&		&		&   &    &    &    &    &    &  & 1 & $\cdots$ \\
	$\cdots$	&		&		&		&   &    &    &    &    &    &  &  &  \\
	total	&	1	&	2	&	4	& 8  &  16  &  32  &  64  &  128  &  256  & 512 & 1 024 & $\cdots$ \\
	

\end{tabular}
\end{center}
\caption{Pascal triangle - Binomial coefficients}
\label{PascalTriangle_BC}

\end{table}


We will use another property of the Diophantine equations generated by functions like $T_3$ and $T_5$.

\begin{thm}
Let the trajectories of the integers (of length L) that are connected to each other by the operations $n/2$ or $(m_in + 1)/2$. The Diophantine equation connecting the first integer $x$ and the last integer $y$ of a sequence can be expressed in the general form $c = by - ax$ where the parameters $a$, $b$ and $c$, always positive, depend on the operations themselves and in which orders they are applied. If $b < a$($\lambda > 1$), $x < y$ and, if $b > a$($\lambda < 1$), $x \ge y$ or $x<y$ when $x$ and $y$ are positive.
	\label{Distribution_x_vs_y}
\end{thm}

\textsl{Proof}

Let $k_1, k_2 = 0, 1, 2, \ldots$ and $k = k_1 + k_2 = L - 1$, with $L \ge 2$.

Then, $a = m_i^{k_2}$, $b = 2^k$ and $c \ge 0$.

As the factors $a$ and $b$ of $x$ and $y$ are prime to one another, the Diophantine equation admits a infinity of solutions to integer values. If $(x_0, y_0)$ is a specific solution, the general solution will be $(x = x_0 + bq, y = y_0 + aq)$, where $q$ is any integer, positive, negative or zero.

Let the equation~\ref{g_k} in the form $y = \lambda x + \rho$, with $\lambda = a/b$ and $\rho$ always positive. A quick examination of this equation allows us to state that if $\lambda > 1$ ($b < a$), $x < y$, and if $\lambda < 1$ ($b > a$), $x \ge y$ or $x < y$. Therefore, two cases are possible, so $b < a$ or $b > a$. Now, let us examine these two cases from the Diophantine equation and its solutions. 


\underline{First case : $b < a$}
 
Of the Diophantine equation $c = by - ax$, as $c$ is always positive and $b < a$, $x$ must always be smaller than $y$ ($x < y$).

\underline{Second case : $b > a$}

Let the general solution

\begin{equation*}
	y = y_0 + aq \phantom{1234} and \phantom{1234} x = x_0 + bq = x_0 + 2^kq,	
\end{equation*}


There are combinations $k = k_1 + k_2$ giving parameters $b$ and $a$ such that $b > a $ ($\lambda < 1$), and maybe integers $x_0 < y_0$ included in the interval $1$ to $2^k$. As $b > a$ and $x_0 < y_0$, beyond a certain value of $q = q_{x > y}$, we will have $x > y$. All integers $x = x_0 + 2^kq$ with $q > q_{x > y}$ will be greater than $y$.

If $x = y$ then $c = x(b - a)$. Since $c$ must always be positive, then $b > a$. $\blacksquare$




For example, for $L = 2$, we have $k = k_1 + k_2 = L - 1 = 1$. Two cases are possible, $k_1 = 1$ and $k_2 = 0$ or, $k_1 = 0$ and $k_2 = 1$. Then, if $m_i = m_3 = 3$, we have $a = 3^{k_2} = 3^0 = 1$ or $a = 3^{k_2} = 3^1 = 3$ and, $b= 2^k = 2$. We write the Diophantine equations

\begin{equation}
	0 = 2y - x \phantom{1234} or \phantom{1234} 1 = 2y - 3x.	
\label{ED_k_eq_2}
\end{equation}

where $(b = 2, a = 1)$ ($b > a$) in the first case and $(b = 2, a = 3)$ ($b < a$) in the other case.

If $(x_0 = 2, y_0 = 1)$ and $(x_0 = 1, y_0 = 2)$ are specific solutions, the general solutions $(x, y)$ are respectively $(2 + 2q, 1 + 1q)$ with ($x > y$), and $(1  + 2q, 2 + 3q)$ with ($x < y$).

This process can be used for any given path length $L$, both for function $T_3$ and $T_5$. 

By examining the Diophantine equation $c = by - ax$ we have deduced that when $b < a$ ($\lambda > 1$), we necessarily have $x < y$ for positive integers.

For the function $T_3$, as $b = 2^k$ and $a = 3^{k_2}$, we can write

\begin{equation*}
	2^k < 3^{k_2} \phantom{1234} or \phantom{1234} k_2 > k \frac{ln2}{ln3}.	
\end{equation*}

For the function $T_5$,

\begin{equation*}
	2^k < 5^{k_2} \phantom{1234} or \phantom{1234} k_2 > k \frac{ln2}{ln5}.	
\end{equation*}

We can summarize these $2$ conditions in only one, that is to say

\begin{equation}
	k_2 > k \theta_i \phantom{1234} where \phantom{1234} \theta_3 = \frac{ln2}{ln3} \phantom{1234} and \phantom{1234} \theta_5 = \frac{ln2}{ln5}.	
\label{condition_x_PP_y}
\end{equation} 
  
          
\section{Distribution function $F(k)$}

Let us define the distribution function $F(k)$ as

\begin{equation}
	F(k) = \lim\limits_{m \rightarrow \infty} \left( 1/m \right) \mu \{ n \le m \phantom{1} | \phantom{1} \chi(n) > k \},
\label{Function_distribution}	
\end{equation} 

where $\mu$ is the number of positive integers $n \le m$ with $m$ that tends towards infinity. $\chi(n)$ is called the "stopping time", and corresponds to the smallest positive integer $k$ such that the iterative application of function $T_{m_i}$ (equation~\ref{mx+1_function}) on a integer $n$ gives the result $T_{m_i}^{(k)}(n) < n$.

We can reformulate the $3x + 1$ conjecture from the definition of stopping time.

$3x+1$ CONJECTURE. Every integer $n \ge 2$ has a finite stopping time.

If we apply $k$ times the function $T_3$ or $T_5$ to each of the positive integers, the distribution function $F(k)$ gives the proportion of the trajectories for which each integer is greater than the starting integer. 

Using the condition obtained at the end of the previous section, it is possible to recursively determine the number of integers which contribute to the distribution function $F(k)$. 

We build the tables~\ref{PascalTriangle_y_PG_x_T_3} and \ref{PascalTriangle_y_PG_x_T_5} (\ref{PascalTriangle_nbr_to_k_EQ_20} and \ref{PascalTriangle_nbr_to_k_EQ_20_5x+1_Problem} in the appendix).

\begin{table}[H]
\begin{center}

\small
\begin{tabular}{c|cccccccccccc}
	
	$k_2 \setminus k$  & 0 & 1 & 2 & 3 & 4 & 5  & 6  & 7  & 8 & 9 & 10 & $\cdots$ \\
		\hline

	0   & 1 & 1\cellcolor[gray]{0.8} & 0 & 0 & 0  & 0  & 0  & 0  & 0  & 0  & 0  & $\cdots$ \\
	1   &   & 1 & 1\cellcolor[gray]{0.8} & 0 & 0  & 0  & 0  & 0  & 0  & 0  & 0  & $\cdots$ \\
	2   &	  &		& 1 & 1 & 1\cellcolor[gray]{0.8}  & 0  & 0  & 0  & 0  & 0  & 0  & $\cdots$ \\
	3	  &		&		&   & 1 & 2  & 2\cellcolor[gray]{0.8}  & 0  & 0  & 0  & 0  & 0  & $\cdots$ \\
	4	  &		&		&		&   & 1  & 3  & 3  & 3\cellcolor[gray]{0.8}  & 0  & 0  & 0  & $\cdots$ \\
	5	  &		&		&		&   &    & 1  & 4  & 7  & 7\cellcolor[gray]{0.8}  & 0  & 0  & $\cdots$ \\
	6	  &		&		&		&   &    &    & 1  & 5  & 12 & 12 & 12\cellcolor[gray]{0.8}  & $\cdots$ \\			
	7	  &		&		&		&   &    &    &    & 1  & 6  & 18 & 30 & $\cdots$ \\
	8	  &		&		&		&   &    &    &    &    & 1  & 7  & 25 & $\cdots$ \\ 
	9	  &		&		&		&   &    &    &    &    &    & 1  & 8  & $\cdots$ \\ 
	10	&		&		&		&   &    &    &    &    &    &    & 1  & $\cdots$ \\ 		
	$\cdots$ 	&		&		&		&   &    &    &    &    &    &    &    & $\cdots$ \\
	$k \theta_3$   &   & 0.63 & 1.26 & 1.89 & 2.52 & 3.15  & 3.79  & 4.41  & 5.05  & 5.68  & 6.31  &  $\cdots$ \\
	\\
	total	&	1	&	1	&	1	& 2  &  3  &  4  &  8  &  13  &  19  & 38 & 64 & $\cdots$ \\
	$2^k$	&	1	&	2	&	4	& 8  &  16  &  32  &  64  &  128  &  256  & 512 & 1 024 & $\cdots$ \\	
	$F(k)$	&	1	&	0.5	&	0.25	& 0.25  &  0.1875  &  0.125  &  0.125  &  0.1016  &  0.0742  & 0.0742 & 0.0625 & $\cdots$ \\	

\end{tabular}
\end{center}
\caption{Pascal triangle - Number of integers $n(i = k_2,j = k)$ by $2^k$ consecutive integers with $\chi > k$ ($a > b$ and $y > x$) - $3x+1$ Problem}
\label{PascalTriangle_y_PG_x_T_3}

\normalsize

\end{table}


\begin{table}[H]
\begin{center}

\footnotesize
\begin{tabular}{c|cccccccccccc}
	
	$k_2 \setminus k$  & 0 & 1 & 2 & 3 & 4 & 5  & 6  & 7  & 8 & 9 & 10 & $\cdots$ \\
		\hline

	0   & 1 & 0 & 0 & 0 & 0  & 0  & 0  & 0  & 0  & 0  & 0  & $\cdots$ \\
	1   &   & 1 & 1 & 0 & 0  & 0  & 0  & 0  & 0  & 0  & 0  & $\cdots$ \\
	2   &	  &		& 1 & 2 & 2  & 0  & 0  & 0  & 0  & 0  & 0  & $\cdots$ \\
	3	  &		&		&   & 1 & 3  & 5  & 5  & 0  & 0  & 0  & 0  & $\cdots$ \\
	4	  &		&		&		&   & 1  & 4  & 9  & 14 & 14 & 14 & 0  & $\cdots$ \\
	5	  &		&		&		&   &    & 1  & 5  & 14 & 28 & 42 & 56 & $\cdots$ \\
	6	  &		&		&		&   &    &    & 1  & 6  & 20 & 48 & 90 & $\cdots$ \\			
	7	  &		&		&		&   &    &    &    & 1  & 7  & 27 & 75 & $\cdots$ \\
	8	  &		&		&		&   &    &    &    &    & 1  & 8  & 35 & $\cdots$ \\ 
	9	  &		&		&		&   &    &    &    &    &    & 1  & 9  & $\cdots$ \\ 
	10	&		&		&		&   &    &    &    &    &    &    & 1  & $\cdots$ \\ 		
	$\cdots$ 	&		&		&		&   &    &    &    &    &    &    &    & $\cdots$ \\

	$k \theta_5$   &   & 0.43 & 0.86 & 1.29 & 1.72 & 2.15  & 2.58  & 3.01  & 3.45  & 3.88  & 4.31  &  $\cdots$ \\
	\\
	total	&	1	&	1	&	2	& 3  &  6  &  10  &  20  &  35  &  70  & 140 & 266 & $\cdots$ \\
	$2^k$	&	1	&	2	&	4	& 8  &  16  &  32  &  64  &  128  &  256  & 512 & 1 024 & $\cdots$ \\
	$F(k)$	&	1	&	0.5	&	0.5	& 0.375  &  0.375  &  0.3125  &  0.3125  &  0.27354  &  0.27354  & 0.27354 & 0.2598 & $\cdots$ \\	

\end{tabular}
\end{center}
\caption{Pascal triangle - Number of integers $n(i = k_2,j = k)$ by $2^k$ consecutive integers with $\chi > k$ ($a > b$ and $y > x$) - $5x+1$ Problem}
\label{PascalTriangle_y_PG_x_T_5}

\normalsize

\end{table}

We used the same numbering as in the  table giving the distribution of dyadic sequences (binomial coefficients) after $k$ iterations, so the rows $i = k_2$ ($0 \le k_2 \le k$) and the columns $j = k$. The values in the gray areas (tables~\ref{PascalTriangle_y_PG_x_T_3} and \ref{PascalTriangle_nbr_to_k_EQ_20}) correspond to the stopping times for the selected $k$ and are not counted in the total.


The number of integers contributing to the distribution function $F(k)$ is therefore obtained recursively from the table which gives the distribution of integers starting the trajectories (different dyadic sequences) by slices of $2^k$ consecutive integers, combined with the condition mentioned at the end of the previous section. In fact, in the sequence of $2^k$ consecutive integers, we retain only those satisfying the condition $k_2 > k \theta_i$ which originates from $a = m_i^{k_2} > b = 2^k$. In the next section we will analyze in detail the evolution of $F(k)$ which respect this condition. We will see that if $k_{2,min}$ is the smallest value of $k_2$ satisfying $k_2 > k \theta_i$ for a given $k$ (for example, $k_{2,min}=4 $ for k=5 in the table~\ref{PascalTriangle_y_PG_x_T_3}), the smallest value of $k_2$ for the following iteration $k + 1$  is $k_{2,min}$ or $k_{2,min}+1$, never more. The evolution of $n(k_2,k)$ is therefore carried out in a staircase, such as $k_{2,min}$ jumping from $0$ (plateau) or $1$. As we will see, the sum of $n(i=k_2, j=k)$ with $\chi>k$ continuously increases with $k$. 

According to the theorem \ref{Distribution_x_vs_y}, we are sure that when $a>b$ we always have $y>x$ meaning that the last integer $y$ of a trajectory resulting from $k$ iterations is greater than the first integer $x$. Since the data $n(i=k_2, j=k)$ from the last tables are calculated recursively (equation~\ref{somme}), one can easily verify that not only the last integer $y>x$, but this will be the case for all the other integers of the trajectory. Then, all the numbers $n(i=k_2, j=k)$ with $k_2$ satisfying the condition $k_2>k\theta_i$ ($a>b$) also satisfy $\chi>k$ and, therefore, contribute to the distribution function $F(k)$.

For all other values of $n(i=k_2, j=k)$, the contribution to $F(k)$ is zero or, at the very least, negligible. According to the theorem \ref{Distribution_x_vs_y}, if $a<b$, then $y\le x$ or $y>x$. The situation $y=x$ corresponds to a loop. Moreover, it has been proved that if the case $y>x$ occurs, it will nevertheless remain limited. If we take as slice of positive consecutive integers the $2^k$ integers of $2$ to $2^k+1$, we know that the dyadic sequences created by each of the integers of this sequence will be different and periodically repeated for all the integers $x + 2^kq$ (generating new slices of $2^k$ consecutive integers), with $q$ any positive integer. If one of the integers $x$ is such that $y>x$ for $a<b$, we will have a finite value $q_{y<x}$ such as $y<x$, and for all $q>q_{y<x}$ we will have $y<x$. Knowing that the analysis of the convergence or not of the trajectories is based on the set of all positive integers, the proportion of the cases that we have just examined ($y>x$ when $a<b$) becomes negligible and thereby their contributions to $F(k)$. This implies that by retaining only the integers satisfying $k_2>k\theta_i$, we have the minimum value of the distribution function. Anyway, even if we took into account the cases $a<b$ with $y>x$ (if non-zero), this would not affect the conclusions that we are going to draw, on the contrary. 


For a given number of iterations $k$, the minimum value of the distribution function is 

\begin{equation}
	F_{min}(k) = \sum_{i > k\theta}^{k} \frac{n(i,k)}{2^k}. 
\label{Function_distribution_new_sum}		
\end{equation} 

Let us take the transition conditions $a_{k_2} > b_{k-1}$ to $a_{k_2} < b_{k}$ for $k$ and $k_2$ given. The first condition corresponds to the trajectories of which all the integers are greater than the starting integer after $k-1$ iterations. The second condition tell us that if we perform another iteration on these integers with same $k_2$ (if the new iteration is an operation on even integers only), the final integer $y$ will be smaller than the starting integer. We then have $\chi = k$, so the number of iterations $k$ satisfying the definition of stopping time (gray areas in the tables~\ref{PascalTriangle_y_PG_x_T_3} and \ref{PascalTriangle_nbr_to_k_EQ_20}). It should be noted that we have verified that for the slice of consecutive positive integers from $2$ to $2^k+1$ the number of integers is exactly the one indicated in the gray areas ($3x+1$ problem) and this, up to with $k=20$ iterations (examples at the end of the appendix).

With all this information, it is easy to build an algorithm giving the different minimum values of the distribution function versus the number of iterations $k$. The first results for the two functions $T_3$ and $T_5$ from $k=10$ to $900$ are presented in the table~\ref{Distribution function F(k)_3(k)_$F_5(k)$}. We can build more detailed tables giving the following $3$ quantities: the number of integers $n(i,k)$ such as $\chi>k$, the number $2^k$ of consecutive positive integers starting from different dyadic sequences and the distribution function $F(k)$ (table~\ref{Distribution function F(k)_3(k)_with_numerator} for the function $T_3$). This type of detailed table brings out the fact that even if the proportion of integers $\chi>k$ gets smaller and smaller, their number tends towards infinity when $k$ tends towards infinity.

\begin{table}
\begin{center}


\begin{tabular}{c|c|c|c|c|c|c}
	\hline
	 $k$ &   $F_3$ & $F_5$ & & $k$ & $F_3$ & $F_5$ \\
	\hline
		10 &  $6.25 \times 10^{-2}$   &  0.25976563   &  & 100 & $2.6396 \times 10^{-4}$ & 0.18060217\\

		20 &  $2.6062 \times 10^{-2}$ &  0.22122192 &  & 200 & $3.0604 \times 10^{-6}$  & 0.17685114 \\	

	  30 &  $1.1894 \times 10^{-2}$ &  0.20572651 &  & 300 & $5.4667 \times 10^{-8}$  & 0.17621811 \\	
		
		40 &  $5.8233 \times 10^{-3}$ &  0.19625785 &  & 400 & $1.1587 \times 10^{-9}$  & 0.17607927 \\
			
		50 &  $3.3167 \times 10^{-3}$ &  0.19116563 &  & 500 & $2.6584 \times 10^{-11}$ & 0.17604048 \\
		
		60 &  $1.9222 \times 10^{-3}$ &  0.18811449 &  & 600 & $6.4455 \times 10^{-13}$ & 0.17603024 \\
		
		70 &  $1.0516 \times 10^{-3}$ &  0.18513014 &  & 700 & $1.5719 \times 10^{-14}$ & 0.17602715 \\
		
		80 &  $6.6440 \times 10^{-4}$ &  0.18317774 &  & 800 & $4.0963 \times 10^{-16}$ & 0.17602622 \\
		
		90 &  $4.1078 \times 10^{-4}$ &  0.18192180 &  & 900 & $1.0837 \times 10^{-17}$ & 0.17602593 \\
	
	\hline
	\end{tabular}
	
	\end{center}
	\caption{Distribution functions $F_3(k)$ and $F_5(k)$}
	\label{Distribution function F(k)_3(k)_$F_5(k)$}
	\end{table}

          
\section{Property of the distribution function $F(k)$}

The inherent properties of the distribution function flow directly from the properties the binomial distribution of integers (Pascal triangles) and the Diophantine equations linking them (via the three theorems). We recall the fact that the column number $j$ corresponds to the exponent $k$ of the parameter $b = 2^k$, so the number of iterations, and the row number $i$ to the exponent $k_2$ of the parameter $a = 3^{k_2}$ or $a = 5^{k_2}$, where $k_2$ is the number of transformations on odd integers. We have $k_2 = 0, 1, \dots, k$.




If $n(i,j) = 0$ for a given combination $i = k_2$ and $j = k$, then $n(i,j) = 0$ for $i = k_2$ fixed and $j > k$. For example, $n(0,1) = 0$ ($a < b$ or $k_2 < k\theta_i$) implies that $n(0, 2)$, $n(0, 3)$, $\cdots$, equal to $0$. The parameter $b = 2^k$ increases (the previous one multiplied by $2$) while the parameter $a = 1$ remains constant, implying that $a$ is always smaller than $b$ and, using the recursion process, $y < x$. Then $\chi < k$ and the new $n(i,j>k) = 0$. In our example, the trajectories are composed of transformations on even integers only (for $k_2=0$), so the first line in the table.  


If $n(i,j) \ne 0$ for a given combination $i = k_2$ and $j = k$, then $n(i,j) \ne 0$ for $i > k_2$ and $j = k$ fixed. For example in the table~\ref{PascalTriangle_y_PG_x_T_3}, $n(6,9) = 12 \ne 0$ ($a > b$ or $k_2 > k\theta_3$) implies that $n(7, 9)$, $n(8, 9)$ and $n(9, 9)$ are different from $0$, because for each new value of $k_2$ the parameter $a$ is the previous one multiplied by $3$. The parameter $a = 3^{k_2} = 3^6$ increases while the parameter $b = 2^k = 2^9$ remains constant, implying that $a$ is always greater than $b$ and $y > x$ (by the theorem~\ref{Distribution_x_vs_y}). Then $\chi > k$ and the new $n(i>k_2,j) \ne 0$.

Now let's look at the possible cases generated by the following two conditions, so $n(i,j) = 0$ and $n(i + 1,j) \ne 0$. Then, $n(i + 1,j + 1) = 0$ or $n(i + 1,j + 1) \ne 0$. The table~\ref{PascalTriangle_nbr_partial} represents examples of these 2 cases.

\begin{table}[H]
\begin{center}
\begin{tabular}{c|ccccccc}
	
	$i \setminus j$ &  & 4 & 5  & 6  & 7  &  & $\cdots$ \\
		\hline

  	   &    &    &     &    &    &    &      \\
	2	   &    & 0 \cellcolor[gray]{0.6} &  0  & 0  & 0  &    & $\cdots$ \\
	3	   &    & 2 \cellcolor[gray]{0.6} &  0 \cellcolor[gray]{0.8}  & 0  & 0  &    & $\cdots$ \\
	4	   &    & 1  &  3 \cellcolor[gray]{0.8} & 3  & 0  &    & $\cdots$ \\
	5	   &    &    &  1  & 4  & 7  &    & $\cdots$ \\
		   &    &    &     &    &    &    &     \\
		   &    &    &     &    &    &    & $\cdots$    \\	

\end{tabular}
\end{center}
\caption{Pascal triangle - Partial view of the table~\ref{PascalTriangle_y_PG_x_T_3} with $\chi$ greater than the number of iterations $k$}
\label{PascalTriangle_nbr_partial}

\end{table}

Indeed, the first condition ($n(i,j) = 0$) implies that $b/a = 2^j/3^i > 1$ and the second condition ($n(i + 1,j) \ne 0$) implies that $b/a = 2^j/3^{i+1} < 1$. By combining these two conditions, we write

\begin{equation*}
	1 < \frac{2^j}{3^i} <3.	
\label{condition_1}
\end{equation*}

Then, the quotient $b/a = 2^{j+1}/3^{i+1}$ for $n(i + 1,j + 1)$ must meet the condition 

\begin{equation*}
	\frac{2}{3} < b/a = \frac{2^j}{3^i} \cdot \frac{2}{3} < 2,	
\label{condition_2}
\end{equation*}

leading to two cases, so $n(i + 1,j + 1) = 0$ if $b/a > 1$ or $n(i + 1,j + 1) \ne 0$ if $b/a < 1$.

The first case represented by the example $n(i=2, j=4) = 0$ and $n(i+1=3, j=4) = 2$ leads to $n(i+1=3, j+1=5) = 0$.


The quotient $b/a = 2^{j+1}/3^{i+2}$ for $n(i + 2,j + 1)$ must meet the condition 

\begin{equation*}
	\frac{2}{3} \cdot \frac{1}{3} < \frac{2^j}{3^i} \cdot \frac{2}{3} \cdot \frac{1}{3} < \frac{2}{3}.	
\label{condition_3}
\end{equation*}

Then, $b < a$ and $n(i + 2,j + 1) \ne 0$, so $n(4,5) = 3$.

The second case represented by the example $n(i=3, j=5) = 0$ and $n(i+1=4, j=5) = 3$ leads to $n(i+1=4, j+1=6) = 3$.

The quotient $b/a = 2^{j+2}/3^{i+1}$ for $n(i + 1,j + 2)$ must meet the condition 

\begin{equation*}
	\frac{4}{3} < \frac{2^j}{3^i} \cdot \frac{2}{3} \cdot \frac{2}{1} < 4.	
\label{condition_4}
\end{equation*}

Then, $b > a$ and $n(i + 1,j + 2) = 0$, so $n(4, 7)= 0$.

Note that we can obtain the previous results by directly analyzing the condition $k_2 > k\theta_i = E + f$ with $E$ an integer ($0 < E < k$) and $f$ a fraction ($0 < f < 1$). 

These properties allow us to follow the evolution of the distribution function $F_{min}(k)$ and de facto, the total number of integers $n$ with $\chi$ greater than the number of iterations $k$, that we call $N_{\chi > k}$. We have proven that if $k_{2,min}$ is the smallest value of $k_2$ satisfying $k_2 > k \theta_i$ ($k_2 = \lceil k \theta_i \rceil$) for a given $k$, the smallest value of $k_2$ for the following iteration $k + 1$  is $k_{2,min}$ or $k_{2,min}+1$, never more. The evolution of $n(k_2,k)$ is therefore carried out in a staircase, such as $k_{2,min}$ jumping from $0$ (plateau) or $1$. As will see, this ensures that even if $F_{min}(k)$ continuously decreases with $k$, the number $N_{\chi > k}$ for each slice of $2^k$ consecutive integers continuously increases, without ever reaching zero.  

We simplify the notation by using $F$ instead $F_{min}$.


Indeed, if the first non-zero value $n(i, k-1)$ after $k-1$ iterations for $i = \lceil (k-1) \theta_i \rceil$ is $a$, the second $b$, the third $c$, $\cdots$, then


\begin{equation}
	F(k-1) = \frac{(a+b+c+\cdots)}{2^{k-1}}.	
\label{dist_1}
\end{equation}

After $k$ iterations we have $n(i,k) = 0$ or $n(i,k) = a$ for the first non-zero value with $i = \lceil k \theta_i \rceil$, and using the property given by the equation~\ref{somme}, the distribution function is


\begin{equation}
	F(k) = \frac{(a+b)+(b+c)+\cdots)}{2^k} = \frac{2(a+b+c+\cdots)-a}{2 \cdot 2^{k-1}} = F(k-1) - \frac{a}{2^k} < F(k-1)
\label{dist_2}
\end{equation}

or,

\begin{equation}
	F(k) = \frac{((a)+(a+b)+(b+c)+\cdots)}{2^k} = \frac{2(a+b+c+\cdots)}{2 \cdot 2^{k-1}} = F(k-1).
\label{dist_3}
\end{equation}

So, $F(k) \le F(k-1)$. The distribution of positive integers $F(k)$, which can simply be called density, decreases constantly with $k$. Nevertheless, the number of integers appearing in the numerator of the function $F(k)$ increases with $k$. If $N_{\chi > k-1} = a + b + c + d + \cdots$ for $k-1$, then $N_{\chi > k} = a + 2(b + c + d + \cdots)$ or $N_{\chi > k} = 2(a + b + c + d + \cdots)$.






Let the number of iterations $k$ approach infinity and evaluate the different limits. Then, the limit of the numerator of $F(k)$, so $N_{\chi>k}(k)$, is

\begin{equation}
	\lim\limits_{k \to \infty}N_{\chi>k}(k)  \to \infty,
\label{limit_num}
\end{equation}

and that of the denominator $2^k$,

\begin{equation}
	\lim\limits_{k \to \infty}2^k  \to \infty.
\label{limit_den}
\end{equation}

For the $3x+1$ problem, we have

\begin{equation}
	\lim\limits_{k \to \infty}F_3(k)  \to 0  \phantom{12} or \phantom{12} \lim\limits_{k \to \infty}1-F_3(k)  \to 1.
\label{limit_3x_plus1}
\end{equation}

For the $5x+1$ problem,

\begin{equation}
	\lim\limits_{k \to \infty}F_5(k)  \to 0.176.  
\label{limit_5x_plus1}
\end{equation}

Even if the values at the numerators and at the denominators of the distribution functions $F_3(k)$ and $F_5(k)$ tend towards infinity, the $\infty / \infty$ ratio is not indeterminate.

The interpretation of these results is presented in the following conclusion.   
\section{Conclusion} 


For the $5x+1$ problem the distribution function $F_5(k)$ tends towards $0.176$ with $k$ tends towards infinity, meaning that more than $17.6$ $\%$ of all positive integers start trajectories such that the first integer $x$ is smaller than any other integer $y$ in the trajectory.

For the $3x+1$ problem, the fact that the distribution function $F_3(k)$ tends towards zero for positive integers when $k$ tends towards infinity, means that the number of these integers such as the stopping time is greater than the number of iterations ($\chi > k$) by slice of $2^k$ consecutive integers tends towards zero. Thus, the average difference between two integers such that $\chi > k$ keeps increasing. However, this does not guarantee that all positive integers, without exception, will have a finite stopping time. It is possible that there remain positive integers which will have an infinite stopping time. Regardless of the value of the number of the iterations $k$ selected, that we can imagine as large as we want, we can say that there will always be trajectories such that the start integer is smaller than any other integer in the trajectory and, the number $N_{\chi>k}(k)$ of these trajectories tends towards infinity as $k$ tends towards infinity.

From the periodicity property that we have updated in this paper, we cannot therefore confirm or deny beyond ant doubt that the $3x+1$ conjecture is true or false. It is perhaps even possible that any non-heuristic approach like the one we has used leads to the conclusion that the conjecture is unprovable. We recall that we used relatively simple mathematical tools (first degree Diophantine equation with two unknowns, binomial distribution and Pascal triangle) and, at no time, we brought into play any probabilistic properties.

\begin{flushright}
\blfootnote{e-mail:roberttremblay02@videotron.ca}
\end{flushright}


\begin{table}
\begin{center}
\small

\rotatebox{90}{
\begin{tabular}{|c|c|c|c|c|c|c}
	\hline
	 $k$ &   number of integers $N_{\chi > k}$ & $2^k$ & $F_3(k)$  \\
	\hline

$\cdots$ &    &  &   \\

		10 &  64                   & $2^{10} = 1 \phantom{1} 024$ & $6.25 \times 10^{-2}$  \\

		20 &  $27 \phantom{1} 328$ & $2^{20} = 1 \phantom{1} 048 \phantom{1} 576$ & $2.6062 \times 10^{-2}$ \\	

	  30 &  $12 \phantom{1} 771 \phantom{1} 274$ & $2^{30} = 1 \phantom{1} 073 \phantom{1} 741 \phantom{1} 824$ & $1.1894 \times 10^{-2}$  \\	
		
		40 &  $6 \phantom{1} 402 \phantom{1} 835 \phantom{1} 000$ & $2^{40} = 1 \phantom{1} 099 \phantom{1} 511 \phantom{1} 627 \phantom{1}  776$ &  $5.8233 \times 10^{-3}$  \\
			
		50 &  $3 \phantom{1} 734 \phantom{1} 259 \phantom{1} 929 \phantom{1} 440$ & $2^{50} = 1 \phantom{1} 125 \phantom{1} 899 \phantom{1}  906 \phantom{1} 842 \phantom{1} 624$ & $3.3167 \times 10^{-3}$  \\
		
		60 &  $2 \phantom{1} 216 \phantom{1} 134 \phantom{1} 944 \phantom{1} 775 \phantom{1} 156$ & $2^{60} = 1 \phantom{1} 152 \phantom{1} 921 \phantom{1}  504 \phantom{1} 606 \phantom{1} 846 \phantom{1} 976$ &  $1.9222 \times 10^{-3}$  \\
		
		70 &  $1 \phantom{1} 241 \phantom{1} 503 \phantom{1} 538 \phantom{1} 986 \phantom{1} 719 \phantom{1} 152$ & $2^{70} = 1 \phantom{1} 180 \phantom{1} 591 \phantom{1}  620 \phantom{1} 717 \phantom{1} 411 \phantom{1} 303 \phantom{1} 424$ &  $1.0516 \times 10^{-3}$  \\
		
		80 &  $803 \phantom{1} 209 \phantom{1} 913 \phantom{1} 882 \phantom{1} 910 \phantom{1} 595 \phantom{1} 105$ & $2^{80} = 1 \phantom{1} 208 \phantom{1} 925 \phantom{1}  819 \phantom{1} 614 \phantom{1} 629 \phantom{1} 174 \phantom{1} 706 \phantom{1} 176$ & $6.6440 \times 10^{-4}$  \\
		
		90 &  $508 \phantom{1} 520 \phantom{1} 069 \phantom{1} 189 \phantom{1} 622 \phantom{1} 659 \phantom{1} 715 \phantom{1} 764$ & $2^{90} = 1 \phantom{1} 237 \phantom{1} 940 \phantom{1}  039 \phantom{1} 285 \phantom{1} 380 \phantom{1} 274 \phantom{1} 899 \phantom{1} 124 \phantom{1} 224$ & $4.1078 \times 10^{-4}$  \\
		
		100 & $302 \phantom{1} 560 \phantom{1} 669 \phantom{1} 500 \phantom{1} 543 \phantom{1} 257 \phantom{1} 546 \phantom{1} 172 \phantom{1} 187$ & $2^{100} = 1 \phantom{1} 267 \phantom{1} 650 \phantom{1} 600 \phantom{1} 228 \phantom{1} 229 \phantom{1} 401 \phantom{1} 496 \phantom{1} 703 \phantom{1} 205 \phantom{1}376$ & $2.3868 \times 10^{-4}$ \\

$\cdots$ &    &  &   \\

200 & $4.9179 \times 10^{54}$ & $2^{200} = 1.6069 \times 10^{60}$ & $3.0604\times 10^{-6}$ \\

300 & $1.1136 \times 10^{83}$ & $2^{300} = 2.0370 \times 10^{90}$ & $5.4667\times 10^{-8}$ \\

400 & $2.9920 \times 10^{111}$ & $2^{400} = 2.5822 \times 10^{120}$ & $1.1587\times 10^{-9}$ \\

500 & $8.7021 \times 10^{139}$ & $2^{500} = 3.2734 \times 10^{150}$ & $2.6584\times 10^{-11}$ \\

600 & $2.6746 \times 10^{168}$ & $2^{600} = 4.1495 \times 10^{180}$ & $6.4455\times 10^{-13}$ \\

700 & $8.2683 \times 10^{196}$ & $2^{700} = 5.2601 \times 10^{210}$ & $1.5719\times 10^{-14}$ \\

800 & $2.7314 \times 10^{225}$ & $2^{800} = 6.6680 \times 10^{241}$ & $4.0963\times 10^{-16}$ \\

900 & $9.1605 \times 10^{253}$ & $2^{900} = 8.4527 \times 10^{270}$ & $1.0837\times 10^{-17}$ \\

$\cdots$ &    &  &   \\
	
	\hline
	\end{tabular}
	}
	\end{center}
	\caption{Number of integers starting trajectories with $\chi>k$ in the interval of $2^k$ consecutive integers, different total number of trajectories ($2^k$) for $k$ iterations, and distribution function $F_3(k)$}
	\label{Distribution function F(k)_3(k)_with_numerator}
	\end{table}

\newpage


\begin{table}[H]
\begin{center}

\scriptsize

\begin{tabular}{c|ccccccccccc}
	
	$k_2 \setminus k$  & $\cdots$ & 11 & 12 & 13 & 14 & 15 & 16 & 17 & 18 & 19 & 20  \\
		\hline
	
	
	0 	& $\cdots$ & 1   & 1   & 1     & 1     & 1     & 1       & 1      & 1      & 1      & 1       \\
	1 	& $\cdots$ & 11  & 12  & 13    & 14    & 15    & 16      & 17     & 18     & 19     & 20      \\
	2 	& $\cdots$ & 55  & 66  & 78    & 91    & 105   & 120     & 136    & 153    & 171    & 190     \\
	3 	& $\cdots$ & 165 & 220 & 286   & 364   & 455   & 560     & 680    & 816    & 969    & 1 140   \\
	4 	& $\cdots$ & 330 & 495 & 715   & 1 001 & 1 365 & 1 820   & 2 380  & 3 060  & 3 876  & 4 845   \\
	5 	& $\cdots$ & 462 & 792 & 1 287 & 2 002 & 3 003 & 4 368   & 6 188  & 8 568  & 11 628 & 15 504  \\
	6 	& $\cdots$ & 462 & 924 & 1 716 & 3 003 & 5 005 & 8 008   & 12 376 & 18 564 & 27 132 & 38 760  \\
	7 	& $\cdots$ & 330 & 792 & 1 716 & 3 432 & 6 435 & 11 440  & 19 448 & 31 824 & 50 388 & 77 520  \\
	8 	& $\cdots$ & 165 & 495 & 1 287 & 3 003 & 6 435 & 12 870  & 24 310 & 43 758 & 75 582 & 125 970 \\
	9 	& $\cdots$ & 55  & 220 & 715   & 2 002 & 5 005 & 11 440  & 24 310 & 48 620 & 92 378 & 167 960 \\
	10	&          & 11  & 66  & 286   & 1 001 & 3 003 & 8 008   & 19 448 & 43 758 & 92 378 & 184 756 \\
	11	&          & 1   & 12  & 78    & 364   & 1 365 & 4 368   & 12 376 & 31 824 & 75 582 & 167 960 \\
	12	&          &		  &  1  & 13    & 91    & 455   & 1 820   & 6 188  & 18 564 & 50 388 & 125 970 \\
	13	&          &		  &     & 1     & 14    & 105   & 560     & 2 380  & 8 568  & 27 132 & 77 520  \\
	14	&          &		  &     &       & 1     & 15    & 120     & 680    & 3 060  & 11 628 & 38 760  \\
	15	&          &		  &     &       &       & 1     & 16      & 136    & 816    & 3 876  & 15 504  \\
	16	&          &		  &     &       &       &       & 1       & 17     & 153    & 969    & 4 845   \\
	17	&          &		  &     &       &       &       &         & 1      & 18     & 171    & 1 140   \\
	18	&          &		  &     &       &       &       &         &        & 1      & 19     & 190     \\
	19	&          &		  &     &       &       &       &         &        &        & 1      & 20      \\
	20	&          &		  &     &       &       &       &         &        &        &        & 1       \\
	
	$\cdots$	&		&		&   &    &    &    &    &    &   &   &   \\

	total	& $\cdots$ &	2 048	&	4 096	& 8 192  &  16 384  &  32 768  &  65 536  &  131 072  &  262 144  & 524 288 & 1 048 576 \\
	

\end{tabular}
\end{center}
\caption{Pascal triangle - Binomial coefficients - $k = 11$ until $k = 20$}
\label{PascalTriangle_BC_to_k_EQ_20}

\end{table}

\normalsize


\begin{table}[H]
\begin{center}

\scriptsize

\begin{tabular}{c|ccccccccccc}
	
	$k_2 \setminus k$  & 11 & 12 & 13 & 14 & 15 & 16 & 17 & 18 & 19 & 20  \\
		\hline
	
	$\cdots$	&		&		&   &    &    &    &    &    &   &   &       \\

	6 	 & 0 & 0 & 0 & 0 & 0 & 0   & 0 & 0 & 0 & 0   \\
	7 	 & 30  & 30\cellcolor[gray]{0.8} & 0 & 0 & 0 & 0  & 0 & 0 & 0 & 0   \\
	8 	 & 55  & 85  & 85\cellcolor[gray]{0.8} & 0 & 0 & 0  & 0 & 0 & 0 & 0  \\
	9 	 & 33  & 88  & 173   & 173 & 173\cellcolor[gray]{0.8} & 0  & 0 & 0 & 0 & 0  \\
	10	 & 9   & 42  & 130   & 303 & 476 & 476\cellcolor[gray]{0.8}   & 0 & 0 & 0 & 0  \\
	11	 & 1   & 10  & 52    & 182   & 485 & 961   & 961 & 961\cellcolor[gray]{0.8} & 0 & 0 \\
	12	 &		  &  1  & 11    & 63    & 245   & 730   & 1 691  & 2 652 & 2 652 & 2 652\cellcolor[gray]{0.8} \\
	13	 &		  &     & 1     & 12    & 75   & 320     & 1 050 & 2 741  & 5 393 & 8 045   \\
	14	 &		  &     &       & 1     & 13    & 88      & 408    & 1 458    & 4 199  & 9 592   \\
	15	 &		  &     &       &       & 1     & 14       & 102     & 510    & 1 968    & 6 167    \\
	16	 &		  &     &       &       &       & 1        & 15      & 117     & 627    & 2 595    \\
	17	 &		  &     &       &       &       &         & 1      & 16      & 133     & 760      \\
	18	 &		  &     &       &       &       &         &        & 1       & 17      & 150       \\
	19	 &		  &     &       &       &       &         &        &        &  1      & 18        \\
	20	 &		  &     &       &       &       &         &        &        &        & 1        \\
	
	$\cdots$	&		&   &    &    &    &    &    &   &   &      \\

	$k \theta_i$    & 6.94 & 7.57 & 8.20 & 8.83 & 9.46  & 10.09  & 10.73  & 11.36  & 11.99  & 12.62  &\\
	\\

	total	&	128	&	226	& 367  &  734  &  1 295  &  2 114  &  4 228  &  7 495  & 14 990 & 27 328  \\
	$2^k$	 &	2 048	&	4 096	&	8 192	& 16 384  &  32 768  &  65 536  &  131 072  &  262 144  &  524 288  & 1 048 576  \\	
	$F(k)$		&	0.0625	&	0.0552	& 0.0448  &  0.0448  &  0.0395  &  0.0323  &  0.0323  &  0.0286  & 0.0286 & 0.0261  \\		


\end{tabular}
\end{center}
\caption{Pascal triangle - Number of integers $n(i = k_2,j = k)$ by $2^k$ consecutive integers with $\chi > k$ ($a > b$ and $y > x$) - $k = 11$ until $k = 20$ - $3x+1$ Problem}
\label{PascalTriangle_nbr_to_k_EQ_20}

\normalsize
\end{table}


\begin{table}[H]
\begin{center}

\scriptsize

\begin{tabular}{c|ccccccccccc}
	
	$k_2 \setminus k$  & 11 & 12 & 13 & 14 & 15 & 16 & 17 & 18 & 19 & 20 \\
		\hline
	
	$\cdots$	&		&		&   &    &    &    &    &    &   &         \\

	5 	 & 56 & 0 & 0 & 0 & 0 & 0 & 0 & 0 & 0 & 0  \\
	6 	 & 146 & 202 & 202 & 0 & 0 & 0   & 0 & 0 & 0 & 0   \\
	7 	 & 165  & 311 & 513 & 715 & 715 & 715  & 0 & 0 & 0 & 0   \\
	8 	 & 110  & 275  & 586 & 1 099 & 1 814 & 2 529  & 3 244 & 3 244 & 0 & 0  \\
	9 	 & 44   & 154  & 429   & 1 015 & 2 114 & 3 928  & 6 457 & 9 701 & 12 945 & 12 945  \\
	10	 & 10   & 54  & 208   & 637 & 1 652 & 3 766   & 7 694 & 14 151 & 23 852 & 36 797  \\
	11	 & 1    & 11  & 65    & 273   & 910 & 2 562   & 6 328 & 14 022 & 28 173 & 52 025  \\
	12	 &      &  1  & 12    & 77    & 350   & 1 260   & 3 822  & 10 150 & 24 172 & 52 345  \\
	13	 &		  &     & 1     & 13    & 90    & 440     & 1 700 & 5 522  & 15 672 & 39 844   \\
	14	 &		  &     &       & 1     & 14    & 104      & 544    &2 244   & 7 766  & 23 438   \\
	15	 &		  &     &       &       & 1     & 15       & 119     & 663    & 2 907    & 10 673    \\
	16	 &		  &     &       &       &       & 1        & 16      & 135     & 798    & 3 705    \\
	17	 &		  &     &       &       &       &          & 1      & 17      & 152     & 950      \\
	18	 &		  &     &       &       &       &          &        & 1       & 18      & 170       \\
	19	 &		  &     &       &       &       &          &        &        &  1      & 19         \\
	20	 &		  &     &       &       &       &          &        &        &        & 1        \\
	
	$\cdots$	&				&   &    &    &    &    &    &   &   &   &    \\

	$k \theta_i$    & 4.73 & 5.17 & 5.60 & 6.03 & 6.46  & 6.89  & 7.32  & 7.75  & 8.18  & 8.61   \\
	\\
	total	& 	532	&	1 008	& 2 016  &  3 830  &  7 660  &  15 320  &  29 925 &  59 850  & 116 456 & 232 912  \\
	$2^k$	& 	2 048	&	4 096	&	8 192	& 16 384  &  32 768  &  65 536  &  131 072  &  262 144  &  524 288  & 1 048 576  \\	
	$F(k)$			&	0.2598	&	0.2461	& 0.2461  &  0.2338  &  0.2338  &  0.2338  &  0.2283  &  0.2283  & 0.2221 & 0.2221  \\		

\end{tabular}
\end{center}
\caption{Pascal triangle - Number of integers $n(i = k_2,j = k)$ by $2^k$ consecutive integers with $\chi > k$ ($a > b$ and $y > x$) - $k = 11$ until $k = 20$ - $5x+1$ Problem}
\label{PascalTriangle_nbr_to_k_EQ_20_5x+1_Problem}

\normalsize
\end{table}

\newpage

EXAMPLES OF RESULTS PREDICTED BY TABLE~\ref{PascalTriangle_y_PG_x_T_3} 


\vspace{2mm}

\underline{Example 1}

Take $k = 3$ and the sequence of consecutive integers $n = 2$ to $2^3 + 1 = 9$.
The trajectories of length $L = k + 1 = 4$ starting with $3$ and $7$  and their dyadic representations, so

\begin{equation*}
(3, 5, 8, 4) \phantom{1,2} and \phantom{1,2} (1, 1, 0), 
\end{equation*}

\begin{equation*}
(7, 11, 17, 26) \phantom{1,2} and \phantom{1,2} (1, 1, 1),
\end{equation*}

are the only ones with a stopping time $\chi > 3$. All the other integers in the sequence, so 2, 4, 5, 6, 8 and 9 produce trajectories with $\chi < 3$. By theorems~\ref{periodicity} and \ref{Distribution_x_vs_y} , we know that all the integers $x + 2^3q$ (where $2 \le x \le 9$) start trajectories with the same finalities. After $3$ iterations, $F(k=3) = 2/8 = 1/4$. So, a quarter of all positive integers $n \ge 2$ have their stopping time greater than $3$ after $3$ iterations, as expected.

\vspace{2mm}

\underline{Example 2}

Take $k = 4$ and the sequence of consecutive integers $n = 2$ to $2^4 + 1 = 17$.

The trajectory starting with $3$ , so

\begin{equation*}
(3, 5, 8, 4, 2) \phantom{1,2} and \phantom{1,2} (1, 1, 0, 0),
\end{equation*}

has a stopping time $\chi =4$.

The trajectories of length $L = k + 1 = 5$ starting with $7$, $11$ and $15$  and their dyadic representations, so

\begin{equation*}
(7, 11, 17, 26, 13) \phantom{1,2} and \phantom{1,2} (1, 1, 1, 0), 
\end{equation*}

\begin{equation*}
(11, 17, 26, 13, 20) \phantom{1,2} and \phantom{1,2} (1, 1, 0, 1), 
\end{equation*}

\begin{equation*}
(15, 23, 35, 53, 80) \phantom{1,2} and \phantom{1,2} (1, 1, 1, 1),
\end{equation*}

as well as $7 + 2^4q$, $\cdots$ are the only ones with a stopping time $\chi > 4$. All the other integers in the sequence, so 2, 4, 5, 6, 8, 9. 10, 12, 13, 14, 16 and 17 produce trajectories with $\chi < 4$. After $4$ iterations, $F(k=4) = 3/16$. So, $3/16$ all positive integers $n \ge 2$ have their stopping time greater than $4$ after $4$ iterations, as expected.

We did the exercice up to $k = 20$ and fully obtained all the result anticipated by tables~\ref{PascalTriangle_y_PG_x_T_3} and \ref{PascalTriangle_nbr_to_k_EQ_20}.

\vspace{4mm}

NOTE

We know that the distribution of the number of integers such of $\chi > k$ for a given number of iterations $k$ is done of $k_2 > (ln2/ln3)k$.

\end{document}